\documentclass[12pt]{article}
\usepackage{amsmath}
\usepackage{amsthm}
\usepackage{amsfonts}

\input boxedeps.tex

\SetepsfEPSFSpecial
\HideDisplacementBoxes

\begin{document}
 
\newtheorem{thm}{Theorem}
\newtheorem{lemma}{Lemma}
\newtheorem{claim}{Claim}
\newtheorem{cor}{Corollary}
\newtheorem{conj}{Conjecture}
\newtheorem{defn}{Definition}

\theoremstyle{remark}
\newtheorem*{rmk}{Remark}

\theoremstyle{plain}
\newtheorem*{Main}{Main Theorem}

\theoremstyle{plain}
\newtheorem*{thmcit}{Theorem}

\title{On Krammer's Representation of the Braid Group}
\author{Matthew G. Zinno\thanks{Partially supported by NSF grants
    DMS-9705019 and DMS-9973232}
  \\Columbia University}
\date{February 11, 2000}
\maketitle

\begin{abstract}
A connection is made between the Krammer representation and the
Birman-Murakami-Wenzl algebra.  Inspired by a dimension argument, a basis is
found for a certain irrep of the algebra, and relations which generate
the matrices are found.  Following a rescaling and change of
parameters, the matrices are found to be identical to those of the
Krammer representation.  The two representations are thus the same,
proving the irreducibility of one and the faithfulness of the other.
\end{abstract}

\section{Introduction}

In [BW], Joan Birman and Hans Wenzl constructed a two-parameter family
of algebras related to braid groups and the Kauffman knot polynomial, 
and analyzed its semisimple
structure.  
The same algebras were discovered simultaneously and independently by
Jun Murakami in [M].
Since the braid group maps homomorphically into this
algebra, the representations on the simple summands of the algebra
give irreducible representations of the braid group.  We will focus on
one particular irreducible representation (for each $n$) out of this
collection, and show that it is identical to an independently
discovered representation.  This is parallel to a similar
identification  which has long
been known between  the Burau representation and a certain summand of
the Hecke algebra.

Daan Krammer, in [K], defined a representation of the braid group
using its interpretation as the automorphism group of the punctured
disk.  Using geometric arguments, he constructed a skein relation
between items called forks.  From the skein relation, he
was able to list a set of algebraic equations describing the action of 
braid group generators on these forks, some of which form a basis for
an invariant module.  These equations thus describe matrix entries for
the representation.  

Using techniques related to the solution of the word problem in the
braid groups in [BKL], Krammer
was then able to prove that this representation was faithful for $n=4$,
thus showing
\begin{thmcit}{\bf (Krammer)}
$B_4$ is linear. \end{thmcit}
This discovery revolutionized the study of braid group linearity, since
for a long time it seemed that the Burau representation was the best
candidate to give a faithful representation of the braid group.


Nevertheless,
Krammer's construction was not widely publicized until Stephen Bigelow
expanded Krammer's result using topological methods:
\begin{thmcit}
{\bf (Bigelow)}
Krammer's representation is faithful for all $n$, thus the braid
groups are linear. \end{thmcit}
This finally answered an important question in braid theory which had
stood since the introduction of the Burau representation in 1935.

The main result of this paper gives a connection between Krammer's
representation and the Birman-Murakami-Wenzl (BMW) algebra:
\begin{Main}
Krammer's representation of the braid group $B_n$ is identical to the
$(n-2) \times 1$ irreducible representation of the BMW 
algebra.
\end{Main}
My proof uses purely algebraic methods, and can
be deduced solely from the equations given in [K] for the action
of the braid group, and the equations given in [BW] for
relations within the algebra.

One immediate consequence of this theorem is
\begin{cor} Krammer's representation is irreducible. \end{cor}
And another follows quickly from Bigelow's result:
\begin{cor} The regular representation of the BMW algebra is faithful.
\end{cor}

Vaughan Jones also discovered this main result, simultaneously and
independently.  However, his methods are different from mine, and
involve a somewhat deeper understanding of the algebraic structure.

The representation which Krammer spelled out
explicitly was originally discovered by Ruth Lawrence in [L] as a
representation of the Hecke algebra, using homology of configuration
spaces similar to those used by Bigelow.  This connection relates to a
topological interpretation of the Jones polynomial which is still in
the process of being explored.  
With the addition of Krammer's presentation of the representation, we
can also make connections to other algebras and knot polynomials.

\subsection{Acknowledgements}

This paper will be part of my Ph.D. thesis.
I'd like to thank my advisor, Joan Birman, 
for her invaluable assistance.  Besides
her everyday commentary and guidance, she brought to my attention many
of the structures and concepts used in this work, and had several
ideas for ways to try to connect them.

I would also like to thank 
Stephen Bigelow, who helped me understand his proof and Lawrence's
contribution, through both correspondence and a visit to Columbia.

I would like to thank Daan Krammer and Vaughan Jones for
correspondence concerning their contributions in this area, and
Justin Roberts, who shared with me his notes and insights from
attending relevant talks given by both 
Bigelow and Jones.

\section{Background and Definitions}

The braid group $B_n$ has many definitions and interpretations.
Perhaps the easiest to see is a pictorial one.  A {\it braid} is a diagram
consisting of two horizontal bars, one at the top and one at the bottom
of the figure, with $n$ nodes on each bar (usually drawn equally
spaced), and $n$ {\it strands}, always running strictly downward, 
connecting the upper and lower nodes.  This figure represents an
isotopy class of embeddings of the strands in
3-space, so the strands are allowed (in fact, required) to cross over
and under each other rather than intersect, and the directions of
these crossings are marked in a conventional way on the diagram.
Among all braids, those which are isotopic in $\mathbb{R} \times I$
are identified.

\begin{figure}[htpb]
\centerline{\BoxedEPSF{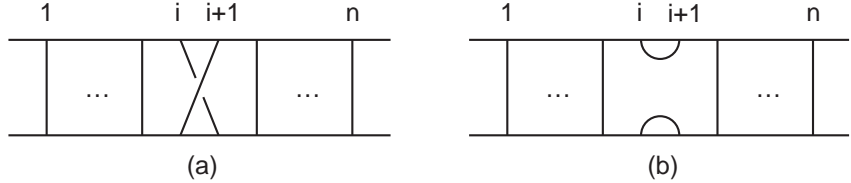 scaled 800}}
\caption {(a) The braid generator $\sigma_i$.
(b) The braid-like algebra element $E_i$.}
\label{fig:generators}
\end{figure}

Braids form a group.  The identity element is the braid with no
crossings, multiplication is concatenation (draw one diagram above the
other, and erase the center bar), and the standard set of generators
(known as {\it Artin generators}) consists of braids $\sigma_i$ 
($1 \leq i < n$), where
the only crossing is that of the $i$th strand under the next strand.  
See Figure \ref{fig:generators}(a).
The relations in the braid group are easily described:
\begin{gather}\label{eqn:braid-rel}
\sigma_i \sigma_j = \sigma_j \sigma_i \text{ for } |i-j|\geq 2 \\
\sigma_i \sigma_{i+1} \sigma_i = \sigma_{i+1} \sigma_i \sigma_{i+1}.
\end{gather}

One of the long-standing questions about the braid group has been
whether it is linear, that is, whether there exists any faithful
representation into a matrix group.  A common method of constructing
representations of the braid group
is to map the braid group homomorphically into a finite-dimensional
algebra, and use the algebra's regular representation.  One such
algebra that was used is the Hecke algebra, a deformation of the
complex symmetric algebra $\mathbb{C}S_n$, and another is the
Birman-Murakami-Wenzl algebra, similarly a deformation of the Brauer algebra.

The BMW algebra $C_n(l,m)$ can be defined on invertible generators
$G_i$ ($1 \leq i <n$), which satisfy the braid relations described
above, and
non-invertible elements 
$E_i$ defined via the formula 
\begin{equation}
G_i + G_i^{-1}=m(1+E_i).
\end{equation}
The additional relations are:
\begin{align}
&E_i E_{i\pm1}E_i = E_i \\
&G_{i\pm1}G_iE_{i\pm1} = E_iG_{i\pm1}G_i = E_iE_{i\pm1} \\
&G_{i\pm1}E_iG_{i\pm1}  = G_i^{-1}E_{i\pm1}G_i^{-1} \\
&G_{i\pm1}E_iE_{i\pm1}  = G_i^{-1}E_{i\pm1} \\
&E_{i\pm1}E_iG_{i\pm1} =E_{i\pm1}G_i^{-1} \\
&G_iE_i=E_iG_i = l^{-1}E_i \\
&E_iG_{i\pm1}E_i = l E_i \\
&E_i^2 = (m^{-1}(l+l^{-1})-1)E_i \\
&G_i^2 = m(G_i+l^{-1}E_i)-1. \label{eqn:last-BW-rel}
\end{align}
See [BW].  
Some of these relations can be deduced from others; however, I am not
concerned here with a minimal presentation.  For those readers who
are, one can be found in [W1].

In addition, it can be deduced that $E_i$ commutes with both $E_j$ and
$G_j$, if $|i-j|\geq 2$.

As in the case of the braid group, we can associate several elements of
the algebra to pictures.  The generators $G_i$ can be identified with
the same braid diagrams (Figure \ref{fig:generators}(a)) 
described above for $\sigma_i$.  The elements
$E_i$ can be identified with similar diagrams where the $i$th node on
the top is joined to the next node on top, and similarly on the bottom
(with all other nodes connected vertically as with $G_i$).  See Figure
\ref{fig:generators}(b). 
Since the
``strands'' do not run strictly downward, this is not a braid
element.  However, we can still ``multiply'' these diagrams by each
other and by braids, to again get isotopy classes of embeddings in
3-space.

The reader may notice that each of the above relations in the algebra
which does not involve addition actually represents an isotopy
relation of these braid-like elements.  Consequently, every monomial
in the algebra can be uniquely represented by such a picture, up to
isotopy.

The irreducible representations of an algebra such as this one
can be identified using
its Bratelli diagram, which encodes the decomposition of each
semisimple algebra $C_n$ into its simple summands.  The regular
representation on each summand (which, as a vector space, is an
invariant subspace) is then an irreducible representation of the
algebra and, by extension, the braid group which maps into it.  
The connections in the diagram between levels encode the 
inductions or restrictions of the representations to the next or
previous level (by adding or ignoring the final generator).

The Bratelli diagram of the Hecke algebra is
the same as the Bratelli diagram of the complex group algebra of the
symmetric group.  So each simple module, thus each irreducible
representation, of $H_n$ is indexed by a Young diagram with $n$ boxes,
and connections from a particular module to each succeeding level
are to the Young diagrams obtained by adjoining one box to the Young
diagram in question.    Define notation as follows:
let $V_{n,\lambda}$
denote the module or representation in the $n$th  level of the Bratelli
diagram, labeled by the Young diagram $\lambda.$


The Hecke algebra is a quotient of the
BMW algebra (obtained by setting $E_i=0$).  Consequently, it
appears as a direct summand of $C_n$, and its
Bratelli diagram is contained in that of the larger algebra.  
The complete structure of the Bratelli diagram of $C_n$, 
as is explained in [W1], is the same as the Brauer algebra
$D_n$:

\begin{thmcit}{\bf (Wenzl)}
(a) $D_n$ is semisimple.

(b) The simple components of $D_n$ are labeled by the set of Young
diagrams with $n-2k$ boxes ($k \in \mathbb{Z}^+$).

(c) If $V_{n,\lambda}$ is a simple $D_n$ module 
it decomposes as a $D_{n-1}$ module into a direct sum of
simple $D_{n-1}$ modules $V_{n-1,\mu}$, where $\mu$ ranges over all
Young diagrams
obtained by removing or (if $\lambda$ contains fewer than $n$ boxes)
adding a box to $\lambda$. 
\end{thmcit}

This means that the Bratelli diagram can be easily constructed using
an inductive method.  The $n=1$ level is a single module, indexed by
the Young diagram consisting of a single box.  The $n=2$ level has
three modules, one indexed by the empty Young diagram, and the other
two indexed by the two Young diagrams of two boxes.  All three are
connected to the module on the previous level.

\begin{figure}[htpb]
\centerline{\BoxedEPSF{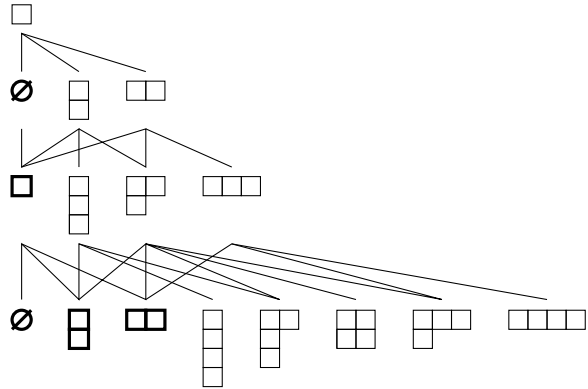 scaled 800}}
\caption {The Bratelli diagram for $C_n$ (up to $n=4$).  
$H_n'$ is in bold.}
\label{fig:Bratelli}
\end{figure}

Starting at level $n=3$, we can construct the levels using
reflections.  Reflecting the $n-2$ level of the diagram, including its
connections to the $n-1$ level, across a line drawn through the $n-1$
level, gives the portion of the $n$ level indexed by Young diagrams
with less than $n$ boxes.  This portion of the algebra, 
following [BW], is given the
notation $H_n'$.  The portion indexed by Young diagrams with $n$ boxes
is then constructed as usual from the modules at the $n-1$ level with
$n-1$ boxes in their Young diagrams. 
This portion of the algebra is
given the notation $H_n$, and is isomorphic to the Hecke algebra.
Consequently, all monomials in the algebra which contain a $E_i$
factor are located in $H_n'$. The Bratelli diagram of $C_n$ up
to $n=4$ is shown in Figure \ref{fig:Bratelli}.  

As with any Bratelli diagram, the dimension of each module (and its
representation) is the sum of the dimensions of the representations it
restricts to, which is thus the number of paths from the top of the
diagram.  For $H_n$, these dimensions can also be found from the hook
length formula.

\section{Dimension Argument}

Since we are claiming that the Krammer representation, which has
dimension $n\choose 2$, is an irreducible representation of $C_n$, it
would be good to know that a representation of the appropriate size
exists:

\begin{thm}In the $n$th level of the 
Bratelli diagram of $C_n$, the representations
labeled by rectangular Young diagrams of shapes $(n-2)\times 1$ and
$1\times(n-2)$ have dimension    $n\choose 2$.
\end{thm}

\begin{proof}
Due to the symmetric and dual nature of the Young diagrams, both
arguments will be similar, and I only need discuss one.  So without loss
of generality, let $\lambda$ denote the $(n-2)\times 1$
Young diagram, consisting of 1 row with $n-2$ boxes.

We can use the Bratelli diagram (see Figure \ref{fig:countingproof}) 
to count the dimension of each of the
representations, as it will be the sum of the dimensions of the
representations on the previous level that this one is connected to.
Since $\lambda$ has fewer than $n$ boxes, this module is part of
$H_n'$, and using the inductive construction of the Bratelli diagram,
we can see that all the connections leading into this representation
come from connections leading out of the representation
$V_{n-2,\lambda} \subset H_{n-2}$ (and, in fact, go to/from the same
representations on the $n-1$ level).  These connections are of two
types: those connecting $V_{n-2,\lambda}$ to modules in $H_{n-1}$, and
those connecting it to modules in $H_{n-1}'$.  The connections to
$H_{n-1}$ 
are easy to see, as they result from the standard Young's
lattice.
There is one representation with a rectangular Young diagram of shape
$(n-1) \times 1$, of dimension 1, and one representation whose Young
diagram has a single box in the second row.  By the hook length
formula, this representation has dimension $n-2$.

\begin{figure}[htpb]
\centerline{\BoxedEPSF{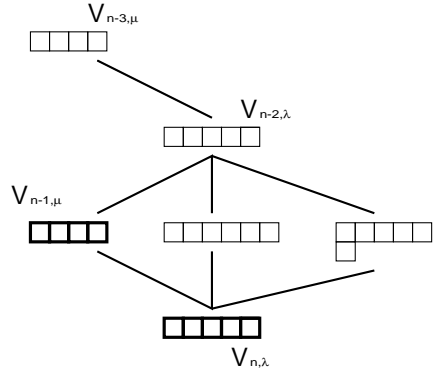 scaled 800}}
\caption {A portion of the Bratelli diagram.}
\label{fig:countingproof}
\end{figure}

Whatever connections exist between $V_{n-2,\lambda}$ and $H_{n-1}'$
again come from reflections of connections from $V_{n-2,\lambda}$ to
$C_{n-3}$.  Of these, there is only one: $V_{n-2,\lambda} \subset
H_{n-2}$ only connects to $H_{n-3}$, and since $\lambda$ has shape
$(n-2) \times 1$, it only connects to the representation labeled by
the Young diagram $\mu$ of shape $(n-3) \times 1$.  So our third and
last downward connection from $V_{n-2,\lambda}$ (and by reflection,
upward connection from $V_{n,\lambda}$) is to $V_{n-1,\mu}$.  Notice
that this is precisely the $n-1$ version of the representation we are
investigating!  Thus by induction, we may assume that this
representation has dimension ${n-1}\choose 2$.  (One may verify that
the claimed formula holds for a sufficiently low base case.)

The dimension of the representation $V_{n,\lambda}$ is therefore
$${{n-1} \choose 2} + (n-2) + 1 ={n\choose 2}.$$
\end{proof}

\section{Explicit form of the representation}

In Jones [J], it is shown that for any braid index $n$, there exists a 
1-dimensional irreducible representation of the Hecke algebra.  More 
specifically, the Hecke algebra has a 1-dimensional invariant subspace 
which is preserved under multiplication by any of the Hecke algebra (or 
braid) generators.  

It is similarly shown in [BW] and [W1] that  $C_n$, 
which contains a subalgebra isomorphic to the Hecke algebra, also contains
a 1-dimensional irreducible representation of 
$C_n$ for any index $n$.  
This is easy to see on the Bratelli diagram as the representation
$V_{n,\lambda}$ (with $\lambda$ as in the proof of Theorem 2).
Thus we again have a 1-dimensional 
invariant subspace which is preserved under multiplication by any of the 
braid generators.  Note that this subspace does not consist of the
same algebra elements as in the 
Hecke algebra.  The quadratic relation is different if we 
consider the full Birman-Murakami-Wenzl algebra, so either subspace is not 
preserved if we use the multiplication rule corresponding to the wrong 
algebra.

Pick an algebra element which is in this 1-dimensional subspace (and 
thus generates it).  We will call it $v$.  Since the actions of the 
braid group on this vector (by left multiplication) result in a 
1-dimensional representation, all such multiplications are scalar 
multiplications by the same factor, which we will call $\kappa$.

To prove that the representation of the braid group $B_n$ given 
explicitly by Krammer is a representation of  $C_n$, 
we will make use of the vector $v$ which corresponds to the braid index 
$n-2$.  Thus, $v \in V_{n-2,\lambda} \subset H_{n-2}$
can be expressed as a linear 
combination of monomials in $C_n$ using only the positive generators
$G_1, \ldots, G_{n-3}$. 
An explicit computation of $v$ in these terms is possible, but not
required for the computations in this paper.  An inductive
construction can be found in [W1].
Importantly, $v$ has the
property that 
$$G_i v = \kappa v \;\;\; \forall i < n-2.$$

\begin{figure}[htpb]
\centerline{\BoxedEPSF{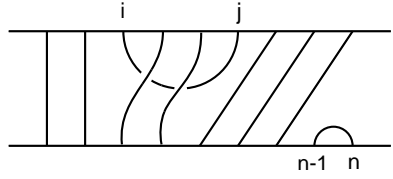 scaled 800}}
\caption {The basis element $T_{ij}$}
\label{fig:T_ij}
\end{figure}

Now we will begin constructing a representation of $C_n$
which we will later show is equivalent to Krammer's 
representation.  A basis for the invariant subspace is the following 
vectors: for $1 \leq i < j \leq n$, let
\begin{equation}\label{eqn:T-ij-def}
\begin{split}
T_{ij} &= \left(\prod_{k=i}^{j-2} G_k\right) \left(\prod_{l=j-
1}^{n-1} E_l\right)v \\
   &= G_i \ldots G_{j-2}E_{j-1}E_j \ldots E_{n-1}v. 
\end{split} \end{equation}
Pictorially, each vector $T_{ij}$ corresponds to the braid-like
diagram where the last two nodes on the bottom bar are connected, and
the $i$th and $j$th nodes on the upper bar are connected (under all
other strands), and the
remaining nodes connect from top to bottom without crossing each
other.  
See Figure \ref{fig:T_ij}.
(Pictures much like these are used in Jones' proof.)

We will occasionally ignore the restriction $i<j$ in order to write more 
general statements.

Now we will see how these vectors behave under a left action by the 
braid group.  Traditionally, Artin braid generators $\sigma_i$ are 
mapped into  $C_n$ under the map 
$\sigma_i \mapsto G_i$; however, that mapping turns out not to work for 
our purposes.  We will instead rescale the braid group, and use the 
mapping 
$\sigma_i \mapsto G_i / \kappa$.  

There are four different types of multiplication to consider (and three 
of them have subtypes which depend on the ordering of the indices):

Type A: $\sigma_i T_{i,i+1}$

Type B: $\sigma_i T_{jk}$, with $\{i,i+1\} \cap \{j,k\}=\emptyset$.

Type C: $\sigma_i T_{i+1,j}$

Type D: $\sigma_i T_{ij}$

For each of these actions, we will map $\sigma_i \mapsto G_i / \kappa$ 
as described above, and expand $T_{ij}$ as defined in equation
(\ref{eqn:T-ij-def}), and simplify 
the resulting expression according to the relations given in equations
(\ref{eqn:braid-rel})-(\ref{eqn:last-BW-rel}).  Most 
of the steps are straightforward; here are lemmas for those that are 
less so:

\begin{lemma}
$E_i G_{i+1} = E_i E_{i+1} G_i^{-1}$
\end{lemma}
\begin{proof}
Both expressions are simplifications of $E_i E_{i+1} E_i G_{i+1}$.
\end{proof}

\begin{lemma}
$E_{i-2} G_i E_{i-1} E_i = E_{i-2} E_{i-1} G_{i-2} E_i$
\end{lemma}
\begin{proof}
Both expressions can simplify to $E_{i-2} G_{i-1}^{-1} E_i$.
\end{proof}

Now to the left action:

A. 
$$\sigma_i T_{i,i+1} = \kappa^{-1} G_i (E_i E_{i+1} \ldots E_{n-1} v)
        =\kappa^{-1}l^{-1} E_i \ldots E_{n-1}v
=\kappa^{-1}l^{-1} T_{i,i+1} $$

B.
Because of the form of $T_{jk}$, the calculation will depend on the order 
of the indices $i,j,k$.

If $i+1<j<k$, then $G_i$ commutes past $T_{jk}$ to multiply by $v$:
\begin{multline}
\sigma_i T_{jk} = \kappa^{-1} G_i (G_j \ldots G_{k-2}E_{k-1}\ldots 
E_{n-1}v) \\
 = \kappa^{-1} (G_j \ldots G_{k-2}E_{k-1}\ldots E_{n-1}) G_i v \\ 
= \kappa^{-1} (G_j \ldots G_{k-2}E_{k-1}\ldots E_{n-1}) (\kappa v) = 
T_{jk}
\end{multline}
 
If $j<i<i+1<k$, then $G_i$ commutes as far as it can, but it is stopped 
when it gets to $G_{i-1}G_i$.  However, the three-term braid relation 
lets it transform into a $G_{i-1}$ and continue commuting to the right:
\begin{multline}
\sigma_i T_{jk} 
= \kappa^{-1} G_i(G_j\ldots G_i \ldots G_{k-2}E_{k-1} \ldots
E_{n-1}v) \\
= \kappa^{-1} G_j\ldots G_{i-2}G_i G_{i-1}G_iG_{i+1} \ldots 
G_{k-2}E_{k-1} \ldots E_{n-1}v \\
= \kappa^{-1} G_j\ldots G_{i-2} G_{i-1}G_i G_{i-1} G_{i+1} \ldots 
G_{k-2}E_{k-1} \ldots E_{n-1}v \\ 
= \kappa^{-1} G_j\ldots G_{i-2} G_{i-1}G_i G_{i+1} \ldots 
G_{k-2}E_{k-1} \ldots E_{n-1} (G_{i-1} v) \\ 
= \kappa^{-1} G_j\ldots G_{k-2}E_{k-1} \ldots E_{n-1}( \kappa v) 
= T_{jk}
\end{multline}

Similarly, if $j<k<i$, then $G_i$ commutes as far as it can, and then we 
use Lemma 2 to commute it the rest of the way, at the expense of 
lowering the index:
\begin{multline}
\sigma_i T_{jk} 
= \kappa^{-1} G_i(G_j \ldots G_{k-2}E_{k-1}\ldots E_i\ldots E_{n-1}v)
\\
= \kappa^{-1} G_j \ldots G_{k-2}E_{k-1}\ldots 
E_{i-2}G_i E_{i-1}E_i\ldots E_{n-1}v \\ 
= \kappa^{-1} G_j \ldots G_{k-2}E_{k-1}\ldots 
E_{i-2} E_{i-1} G_{i-2} E_i\ldots E_{n-1}v \\ 
= \kappa^{-1} G_j \ldots G_{k-2}E_{k-1}\ldots E_{n-1} G_{i-2}v \\ 
= \kappa^{-1} G_j \ldots G_{k-2}E_{k-1}\ldots E_{n-1} (\kappa v ) 
= T_{jk}
\end{multline}

C.  Again, the exact calculations will depend on the order of the 
indices.

If $i+1<j$, the multiplication trivially gives us new indices for $T$.
$$\sigma_i T_{i+1,j}= \kappa^{-1} G_i(G_{i+1} \ldots G_{j-2}E_{j-
1}\ldots E_{n-1}v)= \kappa^{-1}T_{ij} $$

If $j<i$, then the multiplication will behave much as in the latter 
parts of Case B, but with the obstruction this time located at the 
junction of the $G$ terms and the $E$ terms of $T_{j,i+1}$.
\begin{multline}
\sigma_i T_{j,i+1}
=\kappa^{-1} G_i(G_j \ldots G_{i-1}E_i \ldots E_{n-1}v) \\ 
=\kappa^{-1}G_j \ldots G_{i-2}G_i G_{i-1}E_i \ldots E_{n-1}v \\ 
=\kappa^{-1}G_j \ldots G_{i-2}E_{i-1}E_i \ldots E_{n-1}v
=\kappa^{-1} T_{ji} \end{multline}

D.  This time, not only are the calculations dependent on the order of 
the indices, but the results come out differently as well.

If $i+1<j$, then the first term in $T_{ij}$ is $G_i$, so we need to use 
the quadratic relation in  $C_n$.
\begin{multline}
\sigma_i T_{ij}
= \kappa^{-1} G_i(G_i \ldots G_{j-2}E_{j-1}\ldots E_{n-1}v) \\ 
= \kappa^{-1} (m G_i + ml^{-1} E_i -1) 
G_{i+1} \ldots G_{j-2}E_{j-1}\ldots E_{n-1}v \\ 
= \kappa^{-1} (m T_{ij} -T_{i+1,j} + ml^{-1} E_i 
G_{i+1} \ldots G_{j-2}E_{j-1}\ldots E_{n-1}v)
\end{multline}
Two of the three summands we get are recognizable as $T$ vectors, 
but the other begins with an $E$ followed by a product of $G$s.  For 
this, we need to apply Lemma 1.
\begin{multline}
 E_i G_{i+1} G_{i+2}\ldots G_{j-2}E_{j-1}\ldots E_{n-1}v \\
        = E_i E_{i+1} G_i^{-1} G_{i+2}\ldots G_{j-2}E_{j-1}\ldots E_{n-1}v
\\ = E_i E_{i+1} G_{i+2} \ldots G_{j-2}E_{j-1}\ldots E_{n-1}( G_i^{-1} v)
 \\ = E_i E_{i+1} G_{i+2}\ldots G_{j-2}E_{j-1}\ldots E_{n-1}(\kappa^{-1} 
v) 
\end{multline}
It's still not in the standard form for a $T$ vector, but notice what 
one application of Lemma 1 accomplished.  We still have a word made up 
of $E$s followed by $G$s followed by $E$s, and our indices have not 
changed (they're still increasing by one each time).  We've just changed 
the first $G$ that appears into an $E$ of the same index (and gotten a 
$\kappa^{-1}$ scalar out of the multiplication).  Repeated applications 
of Lemma 1 will have the same effect, and we can continue until all our 
$G$ terms have been transformed into $E$ terms.  This gives us the 
vector $E_i \ldots E_{n-1} = T_{i,i+1}$, multiplied by  one $\kappa^{-
1}$ for each $G$ term we had before the first application of the lemma, 
thus $\kappa^{-(j-2-i)}$.  Conclusion:
$$\sigma_i T_{ij} = m \kappa^{-1} T_{ij} - \kappa^{-1} T_{i+1,j} 
+ ml^{-1}\kappa^{i-j+1} T_{i,i+1}$$

If $j<i$, then things come out differently:
\begin{multline}
\sigma_i T_{ji}
= \kappa^{-1} G_i(G_j \ldots G_{i-2}E_{i-1}\ldots E_{n-1}v) \\ 
= \kappa^{-1}G_j \ldots G_{i-2} G_iE_{i-1}E_i\ldots E_{n-1}v
= \kappa^{-1}G_j \ldots G_{i-2} G _{i-1}^{-1}E_i\ldots E_{n-1}v \\ 
= \kappa^{-1}G_j \ldots G_{i-2} (m+mE_{i-1}-G _{i-1})E_i\ldots
E_{n-1}v \\ 
= \kappa^{-1}(m T_{ji}-T_{j,i+1}+
m G_j\ldots G_{i-2}E_i\ldots E_{n-1}v)
\end{multline}
The last of these terms simplifies by commuting each of  the $G$ factors 
past all the $E$ factors.  Each $G$ that reaches $v$ becomes 
multiplication by $\kappa$, and what remains is
$E_i\ldots E_{n-1}v = T_{i,i+1}$.
So we have 
\begin{align*}
\sigma_i T_{ji} &= \kappa^{-1}(m T_{ji}-T_{j,i+1}+m \kappa^{i-1-j} 
T_{i,i+1})  \\ &= m\kappa^{-1} T_{ji}-\kappa^{-1}T_{j,i+1}+
    m \kappa^{i-j-2} T_{i,i+1}
\end{align*}

It is now clear that these vectors form an invariant subspace, of 
dimension $n\choose 2$.  As we have described an action on this space 
by the generators of the braid group, we have a representation of the 
braid group.  This representation has two (complex) parameters, $m$ and 
$l$.  (Recall that $\kappa$ is the eigenvalue of the specific vector 
$v$, so is not a parameter in the same sense.  However, it can be 
adjusted or reset by rescaling, as we did above.)

A final detail is to locate this invariant subspace in $C_n$.

\begin{thm}
The invariant subspace described above is $V_{n,\lambda}$.
\end{thm}
\begin{proof}
From [W2], Prop. (1.2), it is clear that $v \in V_{n-2,\lambda}$ implies 
$T_{n-1,n}=  E_{n-1} v \in
V_{n,\lambda}$.  The cited proposition, applied to the present
situation, states that if $p \in V_{n-2,\lambda}$ is a minimal
idempotent, then $p E_{n-1}$ is a minimal idempotent of  $V_{n,\lambda}$.
For the present purposes, it is not necessary to deal with minimal
idempotents, but notice that since  $V_{n-2,\lambda}$ is
1-dimensional, the minimal idempotent $p$ is a scalar multiple of
$v$.  Therefore, the cited proposition guarantees that a certain scalar
multiple of $v E_{n-1}$ is a minimal idempotent of  $V_{n,\lambda}$, which
is a stronger result than we require.  (Notice also that since $v$ is
written using only the generators $G_1, \ldots, G_{n-3}$, it commutes
with $v$.)  
\end{proof}

\section{The Krammer representation}

We will now recall Krammer's representation from [K], 
written in terms of actions 
on a module with basis $v_{ij}$, with  $1\leq i,j \leq n$ and $i\neq
j$.  Like the 
representation above, it has two (complex) parameters, $q$ and $t$.  
(Stephen Bigelow has found a fascinating interpretation of these 
parameters from a topological perspective.)
Krammer's presentation is longer than that described here, because he 
included formulae for multiplication by band generators (a larger set 
than Artin generators), and because I am taking the liberty to combine 
formulae when the order of the indices does not matter.
\begin{alignat*}{2}
&\sigma_i v_{i,i+1} =tq^2 v_{i,i+1}  & \\
&\sigma_i v_{jk} = v_{jk} &\text{ for } \{i,i+1\} \cap \{j,k\}=\emptyset \\
&\sigma_i v_{i+1,j}=v_{ij} &\text{ for } j\neq i,i+1  \\
&\sigma_i v_{ij} = tq(q-1)v_{i,i+1} + (1-q)v_{ij}+qv_{i+1,j} 
&\text{ if } i+1<j  \\
&\sigma_i v_{ji} = (1-q)v_{ji}+qv_{j,i+1}+q(q-1)v_{i,i+1}
  &\text{ if } j<i 
\end{alignat*}

To show that this is the same as the  representation of  $C_n$
constructed above, we will rescale the basis slightly and set a 
correspondence between parameters:

\begin{thm}
Under the identifications $q=-\kappa^{-2}$, $m=\kappa(1-q)$, 
$l^{-1}=\kappa tq^2$, $v_{ij}=\kappa^{i+j}T_{ij}$, the two actions 
described above are identical.  \end{thm}
\begin{proof}
A.  
$$\sigma_i v_{i,i+1} =\kappa^{2i+1}\sigma_i T_{i,i+1}
=\kappa^{2i}l^{-1} T_{i,i+1} =\kappa^{2i+1} tq^2 T_{i,i+1}
=tq^2 v_{i,i+1}$$

B.
$$\sigma_i v_{jk} = \kappa^{j+k} \sigma_i T_{jk} = \kappa^{j+k}T_{jk} 
= v_{jk}$$

C.
$$\sigma_i v_{i+1,j}=\kappa^{i+j+1} \sigma_i T_{i+1,j}=
\kappa^{i+j} T_{ij} = v_{ij}  $$

D. 
If $i+1<j$, 
\begin{multline}
\sigma_i v_{ij} =\kappa^{i+j}\sigma_i T_{ij}=
 m \kappa^{i+j-1} T_{ij} -\kappa^{i+j-1} T_{i+1,j} 
+ ml^{-1}\kappa^{2i+1)} T_{i,i+1} \\ 
= (1-q)\kappa^{i+j} T_{ij} -\kappa^{i+j-1} T_{i+1,j} 
+ \kappa^2(1-q) tq^2\kappa^{2i+1)} T_{i,i+1} \\ 
= (1-q)v_{ij} -\kappa^{-2}v_{i+1,j} 
+ \kappa^2(1-q) tq^2 v_{i,i+1}
=(1-q)v_{ij}+qv_{i+1,j} + tq(q-1)v_{i,i+1}
\end{multline}

If $j<i$, 
\begin{multline}
\sigma_i v_{ji} =\kappa^{i+j}\sigma_i T_{ji} 
= m\kappa^{i+j-1} T_{ji}-\kappa^{i+j-1}T_{j,i+1}+
  m \kappa^{2i-2} T_{i,i+1} \\ 
= (1-q)\kappa^{i+j} T_{ji}-\kappa^{i+j-1}T_{j,i+1}+
 (1-q) \kappa^{2i-1} T_{i,i+1} \\ 
= (1-q)v_{ji}-\kappa^{-2}v_{j,i+1}+
  \kappa^{-2}(1-q) v_{i,i+1}
= (1-q)v_{ji}+qv_{j,i+1}+q(q-1)v_{i,i+1}
\end{multline}
\end{proof}

\begin{rmk}
The computations in the proof above will actually work under any additional
rescaling of $v_{ij}$ with respect to $T_{ij}$, namely using the
identification 
$v_{ij}=\kappa^{i+j+k}T_{ij}$, for any value of $k$.  For simplicity
in the proof, I set $k=0$, but if these identifications are to be used
for any explicit calculations, I recommend instead using $k=n+1$, so
that the values of the exponent range symmetrically from $-(n-2)$ to
$(n-2)$.  \end{rmk}


\end{document}